\documentclass{amsart}
\usepackage{xcolor}
\usepackage{tikz}
\usetikzlibrary{calc}
\tikzstyle{vertex}=[inner sep = 0pt, minimum width=4pt, fill=black, shape=circle]
\newcommand{\gpoint}[2]{\node[style=vertex, label=#1:$#2$]}
\newcommand{\apoint}[1]{\gpoint{above}{#1}}
\newcommand{\sizeof}[1]{\left\lvert{#1}\right\rvert}

\newtheorem{proposition}{Proposition}[section]

\newtheorem{theorem}[proposition]{Theorem}

\newtheorem{conjecture}[proposition]{Conjecture}

\newcommand{\tauc}{\tau_c}
\newcommand{\nuc}{\nu_c}
\newcommand{\ree}{\mathcal{R}}
\newcommand{\new}[1]{#1}
\title[Packing and Covering Directed Triangles]{Packing and Covering Directed Triangles}
\author{Jessica McDonald \and Gregory J.~Puleo \and Craig Tennenhouse}
\address[Jessica McDonald]{Department of Mathematics and Statistics, Auburn University, Auburn, Alabama, USA 36849, \textup{mcdonald@auburn.edu}}
\address[Gregory J. Puleo]{Department of Mathematics and Statistics, Auburn University, Auburn, Alabama, USA 36849, \textup{gjp0007@auburn.edu}}
\address[Craig Tennenhouse]{Department of Mathematical Sciences, University of New England, Biddeford, ME 04005, \textup{ctennenhouse@une.edu}}
\thanks{The first author is supported in part by NSF grant DMS-1600551}
\begin{document}

\maketitle

\begin{abstract} We prove that if a directed multigraph $D$ has at
  most $t$ pairwise arc-disjoint directed triangles, then there exists
  a set of less than $2t$ arcs in $D$ which meets all directed
  triangles in $D$, except in the trivial case $t=0$. This
  answers affirmatively a question of Tuza from 1990.
\end{abstract}

\section{Introduction}\label{sec:intro}

In the 1980s, Tuza~\cite{TuzaProc, Tuza} posed the following conjecture about packing
and covering triangles in undirected simple graphs (hereafter called graphs). Given a graph $G$,
let $\nu(G)$ be the maximum size of a family of pairwise edge-disjoint triangles in $G$,
and let $\tau(G)$ be the minimum size of an edge set $X$ such that $G-X$ is triangle-free.
Evidently $\tau(G) \geq \nu(G)$, since we are forced to delete at least one edge from
each triangle in a family of edge-disjoint triangles (and these edges must be distinct),
and on the other hand $\tau(G) \leq 3\nu(G)$, since it suffices to delete \emph{all} edges
from each triangle in a maximal family of edge-disjoint triangles. Tuza conjectured
that in fact $\tau(G) \leq 2\nu(G)$ for every graph $G$. As Tuza
observed, this upper bound is sharp if true, and in particular it is achieved by $K_4$ and $K_5$.

The best general result on Tuza's conjecture is due to
Haxell~\cite{Haxell}, who proved that $\tau(G) \leq 2.87\nu(G)$ for
every graph $G$. Other authors have approached the conjecture by
proving that $\tau(G)\leq 2 \nu(G)$ for all graphs in some given
family. Tuza \cite{Tuza} showed that his conjecture holds for all
planar graphs, and Aparna~Lakshmanan, Bujt\'as, and Tuza~\cite{LBT}
showed that it holds for all 4-colorable graphs. The planar result has
been generalized to graphs without $K_{3,3}$-subdivisions (Krivelevich
\cite{Krivelevich}), and then to graphs with maximum average degree
less than $7$ (Puleo \cite{Puleo}). In the case where $G$ is a
$K_4$-free planar graph, the stronger inequality
$\tau(G) \leq \frac{3}{2}\nu(G)$ was proved by Haxell, Kostochka, and
Thomass\'e~\cite{SashaK4}.

Asymptotic, fractional, and multigraph versions of Tuza's conjecture
have also been considered. Yuster \cite{Yuster} proved that
$\tau(G)\leq (2+o(1))\nu(G)$ when $G$ is a dense graph, and this was
shown to be asymptotically tight by Kahn and Baron
\cite{BarKahn}. Yuster \cite{Yuster} also noted that a combination of
results by Krivelevich \cite{Krivelevich} and Haxell and R\"{o}dl
\cite{HaxRod} implies that for any graph $G$ with $n$ vertices,
$\tau(G)< 2 \nu(G)+o(n^2)$. Two fractional versions of Tuza's
Conjecture were proved by Krivelevich~\cite{Krivelevich}. Chapuy,
DeVos, McDonald, Mohar, and Scheide~\cite{CDMMS} tightened one of
these fractional versions, and considered the natural extension of
Tuza's conjecture to multigraphs. Here by multigraph we mean that
multiple edges are permitted, but not loops (they have no
effect on our problem anyways); the definitions of $\mu$ and
$\tau$ are identical to those given in the simple graph case. In
\cite{CDMMS}, planar multigraphs were shown to satisfy Tuza's
conjecture, and $\tau(G)\leq 2.92\nu(G)$ was shown to hold for all
multigraphs $G$.

When posing his conjecture in \cite{Tuza}, Tuza also discussed the problem of packing and covering \emph{directed} triangles. Here by directed multigraph we shall mean any oriented multigraph; by directed graph we shall mean any directed multigraph without parallel arcs in the same direction (but we allow digons, i.e., a pair of arcs $u \to v$ and $v \to u$). Given a directed multigraph $D$, let $\nuc(D)$ denote the maximum size of a family of pairwise arc-disjoint directed triangles, and let $\tauc(D)$
denote the minimum size of an edge set $Y$ such that $D-Y$ has no
directed triangles. Tuza asked: \emph{``Is $\tauc(D)< 2 \nuc(D)$ for every digraph $D$?''}. In this paper we answer this affirmatively with the following theorem.

\begin{theorem}\label{thm:main}
  If $D$ is a directed multigraph with at least one directed triangle, then $\tauc(D) < 2\nuc(D)$.
\end{theorem}

Tuza \cite{Tuza} observed that the rotational 5-tournament $T_5$, pictured in Figure~\ref{fig:5tourn},
satisfies $\tauc(T_5)/ \nuc(T_5)= \tfrac{3}{2}$. Our computational efforts have not yielded any examples with a larger ratio for $\tauc/\nuc$, and in fact we find the following conjecture plausible.

\begin{conjecture}\label{conj:32}
  If $D$ is a directed multigraph, then $\tauc(D) \leq \frac{3}{2}\nuc(D)$.
\end{conjecture}
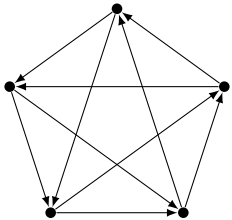
\begin{figure}
  \centering
  \begin{tikzpicture}[rotate=18, scale=1.5]
    \foreach \i in {0,...,4}
    {
      \apoint{} (v\i) at (72*\i : 1cm) {};
    }
    \foreach \i in {0,...,4}
    {
      \draw[-latex] let \n1={int(mod(\i+1,5))} in (v\i) -- (v\n1);
      \draw[-latex] let \n2={int(mod(\i+2,5))} in (v\i) -- (v\n2);
    }
  \end{tikzpicture}
  \caption{The rotational $5$-tournament $T_5$, with $\tauc(T_5) = 3$ and $\nuc(T_5) = 2$.}
  \label{fig:5tourn}
\end{figure}

In \cite{TuzaPerfect}, Tuza proved that if $D$ is a planar oriented graph, then $\tauc(D)=\nuc(D)$. This topic of packing and covering directed triangles appears not to have caught on in the literature however (in contrast to the undirected analogue), and we hope that Conjecture \ref{conj:32} and Theorem \ref{thm:main} may create interest.

\section{Proof of Theorem~\ref{thm:main}}

The main idea of our proof is based on the reducibility argument in
Puleo~\cite{Puleo}. We use induction on $\sizeof{V(D)}$, with trivial
base case when $\sizeof{V(D)} = 1$. \new{Note that in what follows ``triangle'' always means ``directed triangle''.}

  Take any $v \in V(D)$, and define an auxiliary directed multigraph $N$ as
  follows: the vertex set of $N$ is the disjoint union of a set
  $\{s,t\}$ consisting of designated source and sink vertices, as well
  as two sets $W^+$ and $W^-$, where $W^+$ contains a copy $w^+$ of
  each vertex $w \in N^+(v)$, and $W^-$ contains a copy $w^-$ of each
  vertex $w \in N^-(v)$. (Note that if $w \in N^+(v) \cap N^-(v)$, then
  there is a copy of $w$ in \emph{each} of $W^+$ and $W^-$.) Given
  vertices $u^+ \in W^+$ and $z^- \in W^-$, we include the arc
  $u^+ \to z^-$ in $E(N)$ with the same multiplicity as the arc $u \to z$
  in $E(D)$. For each $w^+ \in W^+$, we include the arc $s \to w^+$ in $E(D)$ with
  the same multiplicity as the arc $v \to w$, and for each $w^- \in W^-$, we include the arc $w^- \to t$ in
  $E(N)$ with the same multiplicity as the arc $w \to v$ in $E(D)$.

  Observe that there is a bijection between directed triangles in $D$
  containing $v$, and directed $(s,t)$-paths in $N$; triangle
  $z \to v \to u \to z$ in $D$ corresponds to directed path $su^+z^-t$
  in $N$. Furthermore, two directed triangles in $D$ are arc-disjoint
  if and only if the corresponding paths in $N$ are
  arc-disjoint. (Whenever two triangles use different parallel arcs,
  the corresponding paths have parallel arcs as well.)

  Let $\mathcal{P}$ be a maximum-size set of arc-disjoint
  $(s,t)$-paths in $N$, say with $|\mathcal{P}|=p$.  Let $\ree$ be the
  corresponding set of pairwise arc-disjoint triangles in $D$, all of
  which contain $v$. Each triangle in $\ree$ has exactly one arc that
  is not incident to $v$; let $\ree_v$ be the set consisting of these
  $p$ arcs. 

  Let $X$ be a minimum-size set of arcs in $N$ so that $N-X$ has no
  $(s, t)$-paths. By Menger's Theorem, $|X|=|\mathcal{P}|=p$. Note
  that in $D$, the set $X$ corresponds to a set $X_D$ of $p$ arcs, and
  every triangle incident to $v$ has at least one arc in $X_D$. Let
  $C=X_D\cup \ree_v$, and observe that $C$ is a triangle arc cover of
  every triangle involving $v$ as well as every triangle sharing an
  edge with $\ree$. We have $|C|\leq 2p$, with equality if and only if
  $X_D$ and $\ree_v$ are disjoint.

\new{Let $D'=D-v-\mathcal{R}_v$, and suppose first that} $D'$ has at least one directed triangle. By induction, $\tauc(D')<2\nuc(D')$. Let $\ree'$ be a maximum-size set of edge-disjoint directed triangles in $D'$ and let $C'$ be a minimum-size triangle arc cover in $D'$. By our observations above, note that $C\cup C'$ is a triangle arc cover of $D$, and $\ree' \cup \ree$ is a set of edge-disjoint triangles in $D$. We get that
 \[ |C'\cup C| < 2|\ree'|+2p =2|\ree\cup\ree'|,\]
 as desired.

We may now assume that $D'$ has no directed triangles. In this case,
 $C$ is a triangle arc cover for $D$ with size at most $2p$. Since
 $\ree$ is set of $p$ arc-disjoint triangles in $D$, \new{we may assume that $\tauc(D) = 2p$. We will show that $\nuc(D) \geq p+1$}.

\new{There exists a directed triangle $T_0$ in $D$ that is disjoint from $\ree_v$, since $|\ree_v|=p<\tauc(D)$.  Since $D'$ has no directed triangles, $T_0$ must be incident to $v$; let $e_0$ be the arc of $T_0$ that is not incident to $v$. If $T_0$ has no arcs in common with $\ree$, then $\ree\cup\{T_0\}$ is our desired triangle packing of size $p+1$. Let $\ree^0$ be the set of triangles in $\ree$ with at least one arc in common with $T_0$. Since $e_0\not\in\ree_v$, we know that $|\ree^0|\in\{1, 2\}$. We will show that we can find a set $\mathcal{T}$ of $|\ree^0|$ directed triangles so that $(\ree- \ree^0)\cup\{T_0\}\cup \mathcal{T}$ is a set of $p+1$ arc-disjoint triangles in $D$.}

\new{Let $\ree^0_v$ be the subset of $\ree_v$ that corresponds to $\ree^0$. Consider $D^* = D-v - (\ree_v - \ree^0_v)$. Note that $D^*$ contains at least one triangle, because if not, the arc set  $X_D \cup (\ree_v -\ree^0_v)$ is a triangle arc cover for $D$. Since $D'$ is triangle-free, every triangle in $D^*$ must contain at least one arc from $\ree^0_v$. Hence $|\ree^0_v|\in\{1, 2\}$ implies that $\nu_c(D^*)\in\{1, 2\}$. Let $\mathcal{T}$ be a maximum packing of directed triangles in $D^*$.}

\new{We first claim that $|\mathcal{T}|=|\ree^0|$. If not, then $\nu_c(D^*)=1$ and $|\ree^0|=2$, since every triangle in $\mathcal{T}$ must contain at least one edge from $\ree^0_v$. However $\nuc(D^*) = 1$ implies (by applying the induction hypothesis to $D^*$) that $\tauc(D^*) = 1$, so there is an arc $f^*$ that covers all directed triangles in $D^*$, and hence $X_D \cup (\ree_v - \ree^0_v) \cup \{f^*\}$ is a triangle arc cover in $D$ that is smaller than $C$.}

\new{We now complete our proof by showing that $\mathcal{T}$ is arc-disjoint from $(\ree- \ree^0)\cup\{T_0\}$. Each arc used in this second set of triangles, aside from $e_0$, is either incident to $v$ or from the set $\ree_v-\ree_v^0$. Given that $\mathcal{T}$ is chosen from $D^*$, we need only worry about $e_0$ appearing in some triangle $T \in \mathcal{T}$.
  As observed above, such a $T$ must contain at least one arc from $\ree^0_v$, say the arc $e_1$ from the triangle $R_1 \in \ree^0$. As $R_1$ and $T_0$ share
    an arc incident to $v$, their arcs $e_1$ and $e_0$ either have a common head or a common tail (or both, if they are parallel). Either way, no directed triangle can contain both
    of the arcs $e_1$ and $e_0$, and in particular $T$ cannot contain the arc $e_0$.}

\bibliographystyle{amsplain} \bibliography{dirbib}

\providecommand{\bysame}{\leavevmode\hbox to3em{\hrulefill}\thinspace}
\providecommand{\MR}{\relax\ifhmode\unskip\space\fi MR }
\providecommand{\MRhref}[2]{%
  \href{http://www.ams.org/mathscinet-getitem?mr=#1}{#2}
}
\providecommand{\href}[2]{#2}
\begin{thebibliography}{10}

\bibitem{BarKahn}
Jacob~D. Baron and Jeff Kahn, \emph{Tuza's conjecture is asymptotically tight
  for dense graphs}, Combin. Probab. Comput. \textbf{25} (2016), no.~5,
  645--667. \MR{3531437}

\bibitem{CDMMS}
Guillaume Chapuy, Matt DeVos, Jessica McDonald, Bojan Mohar, and Diego Scheide,
  \emph{Packing triangles in weighted graphs}, SIAM J. Discrete Math.
  \textbf{28} (2014), no.~1, 226--239. \MR{3164555}

\bibitem{Haxell}
P.~E. Haxell, \emph{Packing and covering triangles in graphs}, Discrete Math.
  \textbf{195} (1999), no.~1-3, 251--254. \MR{1663859 (99h:05091)}

\bibitem{HaxRod}
P.~E. Haxell and V.~R\"odl, \emph{Integer and fractional packings in dense
  graphs}, Combinatorica \textbf{21} (2001), no.~1, 13--38. \MR{1805712}

\bibitem{SashaK4}
Penny Haxell, Alexandr Kostochka, and St{\'e}phan Thomass{\'e}, \emph{Packing
  and covering triangles in {$K_4$}-free planar graphs}, Graphs Combin.
  \textbf{28} (2012), no.~5, 653--662. \MR{2964780}

\bibitem{Krivelevich}
Michael Krivelevich, \emph{On a conjecture of {T}uza about packing and covering
  of triangles}, Discrete Math. \textbf{142} (1995), no.~1-3, 281--286.
  \MR{1341453 (96c:05138)}

\bibitem{Puleo}
Gregory~J. Puleo, \emph{Tuza's conjecture for graphs with maximum average
  degree less than 7}, European J. Combin. \textbf{49} (2015), 134--152.
  \MR{3349530}

\bibitem{LBT}
{\relax Aparna Lakshmanan}~S, {\relax Cs}.~Bujt{\'a}s, and {\relax Zs}.~Tuza,
  \emph{Small edge sets meeting all triangles of a graph}, Graphs Combin.
  \textbf{28} (2012), no.~3, 381--392. \MR{2912661}

\bibitem{TuzaProc}
{\relax Zs}olt Tuza, \emph{Finite and infinite sets. {V}ol.\ {I}, {II}},
  Proceedings of the sixth {H}ungarian combinatorial colloquium held in {E}ger,
  {J}uly 6--11, 1981 (Amsterdam) (A.~Hajnal, L.~Lov{\'a}sz, and V.~T. S{\'o}s,
  eds.), Colloquia Mathematica Societatis J\'anos Bolyai, vol.~37,
  North-Holland Publishing Co., 1984, p.~888. \MR{818224 (87a:05005)}

\bibitem{Tuza}
\bysame, \emph{A conjecture on triangles of graphs}, Graphs Combin. \textbf{6}
  (1990), no.~4, 373--380. \MR{1092587 (92j:05108)}

\bibitem{TuzaPerfect}
Zsolt Tuza, \emph{Perfect triangle families}, Bull. London Math. Soc.
  \textbf{26} (1994), no.~4, 321--324. \MR{1302062}

\bibitem{Yuster}
Raphael Yuster, \emph{Dense graphs with a large triangle cover have a large
  triangle packing}, Combin. Probab. Comput. \textbf{21} (2012), no.~6,
  952--962. \MR{2981163}

\end{thebibliography}
\end{document}